\documentclass[smallextended]{svjour3}       
\smartqed  

\usepackage{amsmath,amssymb,xcolor}
\usepackage{graphicx}
\usepackage{bbm}
\usepackage{float}
\usepackage{comment}
\usepackage{enumitem}
\usepackage{authblk}
\usepackage[font=footnotesize,labelfont=bf]{caption}
\usepackage{hyperref}
\usepackage{natbib}
\usepackage{multirow} 
\usepackage{footnote}
\usepackage{xspace}

\newtheorem{condition}{Condition}

\newcommand\independent{\protect\mathpalette{\protect\independenT}{\perp}}
\def\independenT#1#2{\mathrel{\rlap{$#1#2$}\mkern2mu{#1#2}}}

\newcommand{\chd}{CHD\xspace}

\newcommand{\mchd}{MCHD\xspace}

\newcommand{\ohd}{OHD\xspace}

\begin{document}

\title{Bias of the additive hazard model in the presence of causal effect heterogeneity}
\titlerunning{Selection bias of the additive hazard model} 

\author{Richard A.J. Post         \and
 Edwin R. van den Heuvel \and
        Hein Putter 
       
}

\authorrunning{R.A.J. Post, E.R. van den Heuvel, 
and H. Putter} 

\institute{R.A.J. Post \at
              Department of Mathematics and Computer Science, Eindhoven University of Technology, The Netherlands \\
              \email{r.a.j.post@tue.nl}           
        \and 
        E.R. van den Heuvel \at
         Department of Mathematics and Computer Science, Eindhoven University of Technology, The Netherlands
        \and
           H. Putter \at
              Department of Biomedical Data Sciences, Leiden University Medical Center,  The Netherlands
}

\date{Received: date / Accepted: date}

	\maketitle
	
\begin{abstract}Hazard ratios are prone to selection bias, compromising their use as causal estimands. On the other hand, the hazard difference has been shown to remain unaffected by the selection of frailty factors over time. Therefore, observed hazard differences can be used as an unbiased estimator for the causal hazard differences in the absence of confounding. 
However, in the presence of effect (on the hazard) heterogeneity, the hazard difference is also affected by selection. In this work, we formalize how the observed hazard difference (from a randomized controlled trial) evolves by selecting favourable levels of effect modifiers in the exposed group and thus deviates from the causal hazard difference of interest. Such selection may result in a non-linear integrated hazard difference curve even when the individual causal effects are time-invariant. Therefore, a homogeneous time-varying causal additive effect on the hazard can not be distinguished from a constant but heterogeneous causal effect. We illustrate this causal issue by studying the effect of chemotherapy on the survival time of patients suffering from carcinoma of the oropharynx using data from a clinical trial. The hazard difference can thus not be used as an appropriate measure of the causal effect without making untestable assumptions.

    
\end{abstract}

\newpage
\section{Introduction}
Hazard ratios, often obtained by fitting a Cox proportional hazards model \citep{Cox1972}, are the most common effect measures when dealing with time-to-event data. However, the hazard ratio is prone to selection bias due to conditioning on survival and therefore not suitable for causal inference \citep{Hernan2010, Stensrud2017, Aalen2015}. 	It has been recommended to use other, better interpretable, estimands when interested in causal effects \citep{Hernan2010,Bartlett2020,Stensrud2019,Young2020}. Alternatively, using the additive hazard model can avoid interpretation issues \citep{Aalen2015, Martinussen2020}. In the nonparametric model proposed by \citet{Aalen1989}, the hazard rate of time $t$ for individual $i$ in the absence of time-dependent) covariates $\boldsymbol{x}_{i}(t)$ (of dimension $p$) equals
$$\lambda(t \mid \boldsymbol{x}_{i}) = \beta_{0}(t) + \beta_{1}(t) x_{i1}(t) + \hdots + \beta_{p}(t) x_{ip}(t), $$ where the parameters $\beta_{j}(t)$ are arbitrary regression functions, allowing time-varying effects \citep{Aalen2008}. Restricted versions have been proposed by \citet{Lin1994} and \citet{McKeague1994}, where all or some $\beta_{j}(t)$ are assumed to be constant over time. The cumulative regression function, $B_{j}(t)= \int_{0}^{t} \beta_{j}(s) ds$, can reveal whether the effect changes over time, e.g.~in \citet[Example 4.6]{Aalen2008}.   

For cause-effect relations that can be accurately described with Aalen's additive hazard model, the hazard difference is a collapsible measure \citep{Martinussen2013}. Then, in the absence of confounding, the hazard difference is an unbiased estimator for the causal effect, even in the case of unmeasured risk factors \citep{Aalen2015}. For this, it is necessary that the exposure effect on the hazard does not depend on unmeasured individual characteristics (modifiers) and thus are the same for all individuals. However, due to the fundamental problem of causal inference, the effect homogeneity assumption is untestable (even in the absence of any confounding) \citep{Hernan2019}. In our companion paper, we have shown that next to unmeasured risk factors, i.e., frailty,  effect heterogeneity at the level of the individual hazard results in selection bias of observed hazard ratios \citep{Post2022b}. In this work, we study the bias of the hazard difference observed when estimating the causal effect under an additive hazard model but in the presence of heterogeneity of effect (on the hazard). To do so, we formalize the cause-effect relations and define the causal hazard difference. First, we show that in the absence of confounding, the expectation of the  hazard difference observed (in practice) equals a marginalized causal hazard difference. After that, we study how this marginalized causal hazard difference deviates from the causal hazard difference due to selecting individuals with favourable values of the effect modifier. We will illustrate how this selection can result in a non-linear cumulative regression function curve while the individual causal effects are constant. Finally, we discuss the analysis of the effect of treatment on survival with carcinoma of the oropharynx from the clinical trial also studied in \citet{Aalen1989}. 

\section{Notation and hazard differences}
Probability distributions of factual and counterfactual outcomes are defined in the potential outcome framework \citep{Rubin1974, Neyman1990}. Let $T_{i}$ and $A_{i}$ represent the (factual) stochastic outcome and exposure assignment level of individual $i$. Let $T_{i}^{a}$ equal the potential outcome of individual $i$ under an intervention of the exposure to level $a$ (counterfactual when $A_{i} \neq a$). For those more familiar with the do-calculus,  $T^{a}$ is equivalent to $T \mid do(A=a)$ as e.g.~derived in \citet[Equation 40]{Pearl2009} and \citet[Definition 8.6]{Bongers2021}. Throughout this work, we will assume causal consistency, i.e.~if $A_{i}=a$, then $T_{i}^{a}=T_{i}^{A_{i}}=T_{i}$. Causal consistency implies that potential outcomes are independent of the exposure levels of other individuals (no interference). The hazard rate of the potential outcome can vary among individuals due to heterogeneity in risk factors $U_{0}$ as also considered by \citet{Aalen2015}. The hazard difference of the potential outcomes with and without ($a=0)$ exposure might also vary among individuals due to an effect modifier $U_{1}$. Therefore, the hazard rate of individual $i$ at time $t$ of the potential outcome under exposure to level $a$ is a function of $U_{0i}$ and $U_{1i}$ and thus random and equals
\begin{equation}\label{eq:chrate}
\lambda_{i}^{a}(t)=\lim_{h\rightarrow 0}h^{-1}\mathbb{P}\left(T_{i}^{a} \in [t,t+h) \mid T_{i}^{a} \geq t, U_{0i}, U_{1i}\right).\end{equation} The hazard of the potential outcome $T_{i}^{a}$ can be parameterized with a function that depends on $U_{0i}$, $U_{1i}$ and $a$. We focus on cause-effect relations that can be parameterized with a structural causal model (SCM), 
that consists of a joint probability distribution of $(N_{A}, U_{0}, U_{1}, N_{T})$ and a collection of structural assignments (for more details see \citet[Section 2]{Post2022b}) such that

\begin{center}
	\fbox{%
		\parbox{0.9\linewidth}{%
			\begin{align}\label{CH6SCMsurv}
			A_{i}&:=f_{A}(N_{Ai})   \\\nonumber
			\lambda_{i}^{a}(t)&:= f_{0}(t,U_{0i})+f_{1}(t,U_{1i},a) \\ \nonumber
			T_{i}^{a}&:= \min \{t>0:e^{-\Lambda_{i}^{a}(t)}\leq N_{Ti}\},
			\end{align} 			where $\Lambda_{i}^{a}(t)=\int_{0}^{t}\lambda_{i}^{a}(s)ds$, $f_{1}(t,U_{1i},0)=0 $,  $N_{Ai},N_{Ti}\sim \text{Uni}[0,1]$, while $ U_{0i}, U_{1i} \independent N_{Ti}$ and $f_{A}$ is the inverse cumulative distribution function of $A$.}}
\end{center}

\noindent If $\lambda_{i}^{a}(t)-\lambda_{i}^{0}(t) = f_{1}(t,a)$, then the effect is equal for each individual, i.e.~effect homogeneity, and Aalen's additive hazard model is well-specified. Otherwise, $\lambda_{i}^{a}(t)-\lambda_{i}^{0}(t)$  differs among individuals so that $\mathbb{E}[\lambda_{i}^{a}(t)-\lambda_{i}^{0}(t)]$ will be the estimand of interest. The latter contrast equals the difference between the expected hazard rate in the world where everyone is exposed to $a$ and the world without exposure. It will be referred to as the causal hazard difference (\chd) defined in Definition \ref{CH6def:CHD}. 

\begin{definition}{\textbf{Causal hazard difference}}\label{CH6def:CHD}
	The \chd for cause-effect relations that can be parameterized with SCM \ref{CH6SCMsurv} equals
	\begin{align*}
 \mathbb{E}\left[\lambda_{i}^{a}(t)\right]-\mathbb{E}\left[\lambda_{i}^{0}(t)\right]= \mathbb{E}[f_{1}(t,U_{1},a)]=&
    \int\lim_{h\rightarrow 0}h^{-1}\mathbb{P}\left(T^{a} \in [t,t+h) \mid T^{a} \geq t, U_{0}, U_{1} \right)dF_{U_{0}, U_{1}}\\
   &  -\int\lim_{h\rightarrow 0}h^{-1}\mathbb{P}\left(T^{0} \in [t,t+h) \mid  T^{0}\geq t, U_{0}  \right)dF_{U_{0}}.
	\end{align*} 
\end{definition} \noindent Throughout this work, we abbreviate the Lebesque-Stieltjes integral of a function $g$ with respect to probability law $F_{X}$, $\int g(x) dF_{X}(x)$, as $\int g(X) dF_{X}$. 

The \chd thus equals the difference of hazard rates of the potential outcomes marginalized over the population distribution of $(U_{0}, U_{1})$. 
However, the distribution of $(U_{0},U_{1})$ among survivors will differ over time (in all worlds), i.e.~$(U_{0},U_{1}) \overset{d}{\neq} (U_{0},U_{1})\mid T^{a} \geq t$. In turn, $(U_{0},U_{1})\mid T^{a} \geq t$ and $(U_{0},U_{1})\mid T^{0} \geq t$ can differ in distribution. Therefore, the hazard rates, and thus the hazard difference, for the factual outcomes are affected by these conditional distributions of $U_{0}$ and $U_{1}$. We refer to the expected value of the difference of the observed hazards as the observed hazard difference (\ohd) presented in Definition \ref{CH6def:OHD}. 

\begin{definition}{\textbf{Observed hazard difference}}\label{CH6def:OHD}
	The \ohd at time $t$ equals
	\begin{align*}
	&\lim_{h\rightarrow 0}h^{-1}\mathbb{P}\left(T \in [t,t+h) \mid T \geq t, A=a \right)-\lim_{h\rightarrow 0}h^{-1}\mathbb{P}\left(T \in [t,t+h) \mid T \geq t, A=0 \right)\\ =& 
	\int\lim_{h\rightarrow 0}h^{-1}\mathbb{P}\left(T \in [t,t+h) \mid T \geq t, A=a, U_{0}, U_{1} \right)dF_{U_{0}, U_{1}\mid T \geq t, A=a}\\
	&-\int\lim_{h\rightarrow 0}h^{-1}\mathbb{P}\left(T \in [t,t+h) \mid  T \geq t, A=0, U_{0}  \right)dF_{U_{0}\mid T \geq t, A=a}.
	\end{align*}
\end{definition}

 \noindent For a randomized controlled trial (RCT), without informative censoring,  the \ohd can be easily expressed in terms of potential outcomes. \noindent By causal consistency \begin{equation}\mathbb{P}\left(T \in [t,t+h) \mid T \geq t, A=a  \right)=\mathbb{P}\left(T^{a} \in [t,t+h) \mid T^{a} \geq t, A=a  \right).  \end{equation} \noindent By design of the trial $A \independent T^{a}$ (in SCM \eqref{CH6SCMsurv} equivalent to $N_{A} \independent U_{0}, U_{1}, N_{T})$, so \begin{equation}\mathbb{P}\left(T^{a} \in [t,t+h) \mid T^{a} \geq t, A=a  \right)=\mathbb{P}\left(T^{a} \in [t,t+h) \mid T^{a} \geq t \right).  \end{equation} The \ohd at time $t$ is thus equal to 
 \begin{equation}\label{CH6:MCHD} \mathbb{P}\left(T^{a} \in [t,t+h) \mid T^{a} \geq t \right) - \mathbb{P}\left(T^{0} \in [t,t+h) \mid T^{0} \geq t \right). \end{equation}  We refer to this effect measure as the marginal causal hazard difference (\mchd) defined in Definition \ref{CH6def:MCHD}. 
 
 \begin{definition}{\textbf{Marginal causal hazard difference}}\label{CH6def:MCHD}
	The \mchd at time $t$ for cause-effect relations that can be parameterized with SCM \ref{CH6SCMsurv} equals
	\begin{align*}
	&\lim_{h\rightarrow 0}h^{-1}\mathbb{P}\left(T^{a} \in [t,t+h) \mid T^{a} \geq t  \right)-\lim_{h\rightarrow 0}h^{-1}\mathbb{P}\left(T^{0} \in [t,t+h) \mid T^{0} \geq t \right)\\ =& 
	\int\lim_{h\rightarrow 0}h^{-1}\mathbb{P}\left(T^{a} \in [t,t+h) \mid T^{a} \geq t, U_{0}, U_{1} \right)dF_{U_{0}, U_{1}\mid T^{a} \geq t}\\
	&-\int\lim_{h\rightarrow 0}h^{-1}\mathbb{P}\left(T^{0} \in [t,t+h) \mid  T^{0}\geq t, U_{0}  \right)dF_{U_{0}\mid T^{0} \geq t}.
	\end{align*}
\end{definition} \noindent As the integration in Definition \ref{CH6def:CHD} is with respect to the population distribution of $(U_{0}, U_{1})$, the \ohd can deviate from the \chd. However, it has been explained by \citet{Aalen2015} that $U_{0} \independent A \mid T \geq t$ so that for degenerate $U_{1}$, the \mchd equals the \chd so that the latter can be unbiasedly estimated from RCT data. In this paper, we formalize the \mchd (and thus, in the absence of confounding, the \ohd) in case of effect heterogeneity (non-degenerate $U_{1}$). 

\section{Results}\label{CH6sec:hazdiff}
In the remainder of the paper, we will focus on binary exposures such that $a \in \{0, 1\}$. It is known from \citet{Aalen2015} that when $$\lambda_{i}^{a}(t) = f_{0}(t,U_{0i})+f_{1}(t,a),$$ $U_{0i} \independent A_{i} \mid T_{i} \geq t$. As a result, the selection of frailty factors $U_{0}$ over time is similar for exposed and unexposed individuals and does not result in bias. In the absence of confounding ($T^{a} \independent A$), the independence and causal consistency imply
\begin{align}
U_{0}\mid \{T \geq t\} &\overset{d}{=} U_{0} \mid T \geq t, A=a\\ \nonumber
& \overset{d}{=} U_{0} \mid T^{a} \geq t, A=a\\ \nonumber
& \overset{d}{=} U_{0} \mid T^{a} \geq t.
\end{align} Thus \citet{Aalen2015}, implicitly  showed that $U_{0} \mid T^{1} \geq t \overset{d}{=} U_{0} \mid T^{0} \geq t$. In Lemma \ref{CH6l5.1}, we show that the additive frailty remains exchangeable in the presence of effect heterogeneity at the hazard scale. 
\begin{lemma}\label{CH6l5.1}
	If the cause-effect relations of interest can be parameterized with SCM \eqref{CH6SCMsurv}, where 
	\begin{equation*}\lambda_{i}^{a}(t):= f_{0}(t,U_{0i})+f_{1}(t,U_{1i},a),
	\end{equation*} and $U_{0i} \independent U_{1i}$ then, 
	\begin{equation*}
	\mathbb{E}\left[f_{0}(t,U_{0})\mid T^{1} \geq t\right]=	\mathbb{E}\left[f_{0}(t,U_{0})\mid T^{0} \geq t\right].
	\end{equation*}
\end{lemma} \noindent Note that while $\mathbb{E}\left[f_{0}(t,U_{0})\mid T^{1} \geq t\right] = \mathbb{E}\left[f_{0}(t,U_{0})\mid T^{0} \geq t\right]$, $\mathbb{E}\left[f_{0}(t,U_{0})\mid T^{a} \geq t\right] \neq \mathbb{E}\left[f_{0}(t, U_{0})\right]$ as the conditional expectations will decrease over time (survival of less susceptible individuals). The selection of levels of $U_{0}$ is thus similar in the exposed and unexposed worlds. It indeed follows from Lemma \ref{CH6l5.1} that in the absence of effect heterogeneity, the \ohd from an RCT is an unbiased estimator of the \chd as shown by \citet{Aalen2015}. However, if heterogeneity exists, there will also be a selection of $U_{1}$ in the exposed world where individuals with more favourable levels of $U_{1}$ are more likely to survive. As a result of this selection, the \ohd over time no longer represents the (population) average effect. In this work, we consider hazard functions that satisfy Condition \ref{CH6regconda}.

\begin{condition}{\textbf{Hazard without infinite discontinuity}}\label{CH6regconda}
	$$\exists \tilde{{}h}>0 \text{ such that }\forall h^{*} \in (0,\tilde{{}h}):  \mathbb{E}\left[f_{0}(t+h^{*},U_{0})+f_{1}(t+h^{*},U_{1},a)\mid T^{a}\geq t\right]<\infty$$
\end{condition}

\noindent In Theorem \ref{CH6th52}, we show that the \ohd equals $\mathbb{E}\left[f_{1}(t,U_{1},1)\mid T^{1} \geq t\right]$ in absence of confounding for hazard functions that satisfy Condition \ref{CH6regconda}, and thus deviates from the \chd equal to $\mathbb{E}[f(_{1}(t,U_{1},1)]$.  

\begin{figure}[H]
		\centering
		\resizebox{0.65\textwidth}{!}{\includegraphics[width=\textwidth]{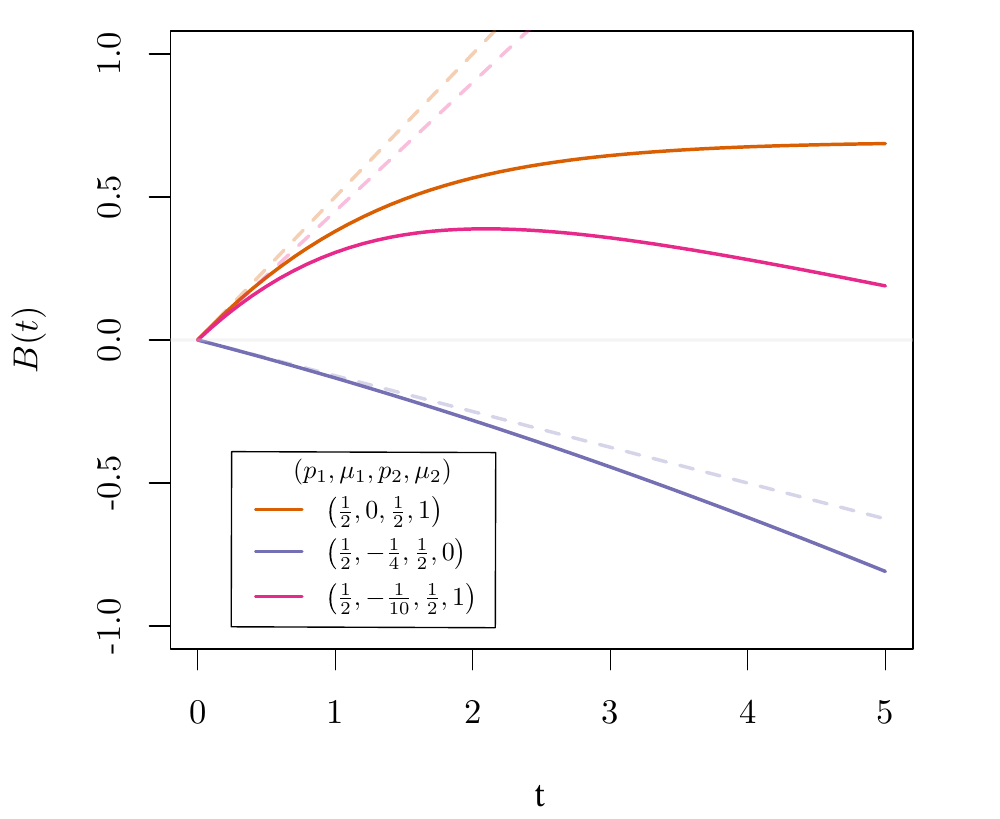}} %
		\caption{$\int_{0}^{t} \mathbb{E}\left[U_{1} \mid T^{1}\geq s\right] ds$, when $U_{1}\sim \text{BHN}(p_{1}, \mu_{1}, p_{2}, \mu_{2})$ (solid), and $g(t)=t \mathbb{E}[U_{1}]$ (dashed). }\label{CH6FigBHN}	
\end{figure}

\begin{theorem}\label{CH6th52}
	If the cause-effect relations of interest can be parameterized with SCM \eqref{CH6SCMsurv}, where 
	\begin{equation*}\lambda_{i}^{a}(t):= f_{0}(t,U_{0i})+f_{1}(t,U_{1i},a),
	\end{equation*} $A_{i} \independent T_{i}^{a}$ and Condition \ref{CH6regconda} holds
	then
	\begin{align*}
	&\lim_{h\rightarrow 0}h^{-1}\mathbb{P}\left(T \in [t,t+h) \mid T \geq t, A=1 \right)-\lim_{h\rightarrow 0}h^{-1}\mathbb{P}\left(T\in [t,t+h) \mid T \geq t, A=0 \right)\\
	=&\mathbb{E}[ f_{1}(t,U_{1},1)\mid T^{1} \geq t] + \mathbb{E}[ f_{0}(t,U_{0})\mid T^{1} \geq t] -\mathbb{E}[f_{0}(t,U_{0}) \mid T^{0} \geq t],
	\end{align*}
	which equals $\mathbb{E}[ f_{1}(t,U_{1},1)\mid T^{1} \geq t]$ when $U_{0} \independent U_{1}$. 
\end{theorem} \noindent In case $f_{1}(t,U_{1i},a) = U_{1i}f_{1}(t,a)$, $\mathbb{E}\left[f_{1}(t,U_{1},1)\mid T^{1} \geq t\right]=f_{1}(t,a)\mathbb{E}[U_{1}\mid T^{1} \geq t]$. By Definition \ref{CH6def:CHD}, the \chd equals $f_{1}(t,a)\mathbb{E}[U_{1}]$. So, the difference with the \ohd varies over time and equals. \begin{equation}f_{1}(t,a)\left(\mathbb{E}[U_{1}]-\mathbb{E}[U_{1}\mid T^{1} \geq t]\right).\end{equation} For this case, the conditional expectation $\mathbb{E}[U_{1}\mid T^{1} \geq t]$ can be expressed in terms of the Laplace transform of $U_{1}$, as shown in Lemma \ref{CH6lemLP}. 
\begin{lemma}\label{CH6lemLP}
	If the cause-effect relations of interest can be parameterized with SCM \eqref{CH6SCMsurv}, where 
	\begin{equation*}f_{1}(t,U_{1i},1) = U_{1i}f_{1}(1,t), 
	\end{equation*} and $U_{0i} \independent U_{1i}$, then 
	\begin{equation}
	\mathbb{E}\left[U_{1} \mid T^{1} \geq t\right]   
	=-\frac{\mathcal{L}_{U_{1}}^{'}(\int_{0}^{t}f_{1}(1,s)ds)}{\mathcal{L}_{U_{1}}(\int_{0}^{t}f_{1}(1,s)ds)},
	\end{equation} where $\mathcal{L}_{U_{1}}(c)=\mathbb{E}\left[\exp(-cU_{1})\right]$ with derivative $\mathcal{L}_{U_{1}}^{'}(c)$.
\end{lemma} \noindent We will apply this lemma to illustrate how effect heterogeneity can affect the cumulative regression function curve in case the causal effect is constant for each individual ($f_{1}(t,U_{1i},1)=U_{1i}$). Let the additive hazard effect modifier $U_{1}$ equal $\mu_{1}$ ($\leq 0$, for individuals that benefit) with probability $p_{1}$, $\mu_{2}$ ($\geq 0$, for individuals that are harmed) with probability $p_{2}$ or $0$ (for individuals that are not affected). We denote this distribution as the Benefit-Harm-Neutral, $\text{BHN}(p_{1},\mu_{1},p_{2},\mu_{2})$, distribution. Note that it is necessary that $\forall t: \mathbb{P}(f_{0}(t, U_{0})<\mu_{1})=0$ to guarantee that the hazard rate is positive for each individual. The Laplace transform of $U_{1}$ equals, 
$$p_1 e^{-c \mu _1}+p_2 e^{-c \mu _2}+\left(1-p_1-p_2\right).$$ So, by Lemma \ref{CH6lemLP} with $f_{1}(1,t)=1$, the \ohd from an RCT is equal to
$$ \mathbb{E}\left[U_{1} \mid T^{1} \geq t \right]=\frac{\mu _1 p_1 e^{-t \mu _1}+\mu _2 p_2 e^{-t \mu _2} }{p_1 e^{-t \mu _1}+p_2 e^{-t \mu _2}+\left(1-p_1-p_2\right) }.$$ The \ohd then deviates from the constant \chd equal to $\mathbb{E}[U_{1}]=p_{1}\mu_{1} + p_{2}\mu_{2}$. For cause-effect systems that meet the conditions of Theorem \ref{CH6th52}, the expectation of the integrated \ohd (and thus of the integrated \mchd) equals 
$$ B(t) = \int_{0}^{t} \mathbb{E}\left[U_{1} \mid T^{1}\geq s\right] ds =  -\log\left(p_{1}(\exp(-t\mu_{1})-1)+p_{2}(\exp(-t\mu_{2})-1)+1\right).$$ Even though the \chd is constant, $B(t)$ will not be linear and deviate from the function $g(t)=t \mathbb{E}[U_{1}]$ due to the selection effect (of $U_{1}$) over time. Three types of curves could be observed as shown in Figure \ref{CH6FigBHN} where for illustration, $p_{1}=p_{2}=0.5$, so there exist only two levels for the modifier $U_{1}$. 

First of all, let's consider the case where the exposure harms some individuals (for which $U_{1i}=1$) while others don't respond to the exposure at all ($U_{1i}=0$), see the orange line in Figure \ref{CH6FigBHN}. Initially, $B(t)$ evolves as $t\mathbb{E}[U_{1}] = 0.5 t$. 
However, the individuals harmed by the exposure are less likely to survive over time, so the curve's derivative decreases. In the end, only individuals with $U_{1i}=0$ are expected to survive so that $B(t)$ remains constant. Concluding that the exposure initially harms but loses effect over time is false for this case as the effect remains constant for each individual. 

Secondly, when some individuals do benefit $(U_{1i}=-0.25)$ while others are not affected ($U_{1i}=0$), the derivative of $B(t)$ evolves from $-0.125$ to $-0.25$ at the moment only those that benefit are expected to survive, as illustrated with the purple line in Figure \ref{CH6FigBHN}. But, again, the effect for an individual is constant and does not become more beneficial over time. 

Finally, and probably most interesting, different individuals in the population might have opposite effects ($U_{1i}=1$ or $U_{1i}=-0.1$), as illustrated with the pink line in Figure \ref{CH6FigBHN}. Initially, the integrated hazard differences increase as the expected effect is harmful. However, over time those individuals with $U_{1i}=-0.1$ are more likely to survive so that $\mathbb{E}[U_{1}\mid T^{1}\geq 1]$ changes sign. Finally, only those that benefit are expected to survive, and the curve decreases with a derivative equal to $-0.1$. For this example, it would be false to conclude that the exposure first harms but becomes beneficial over time. 

 \noindent

\begin{figure}[H]
		\centering
		\resizebox{0.65\textwidth}{!}{\includegraphics[width=\textwidth]{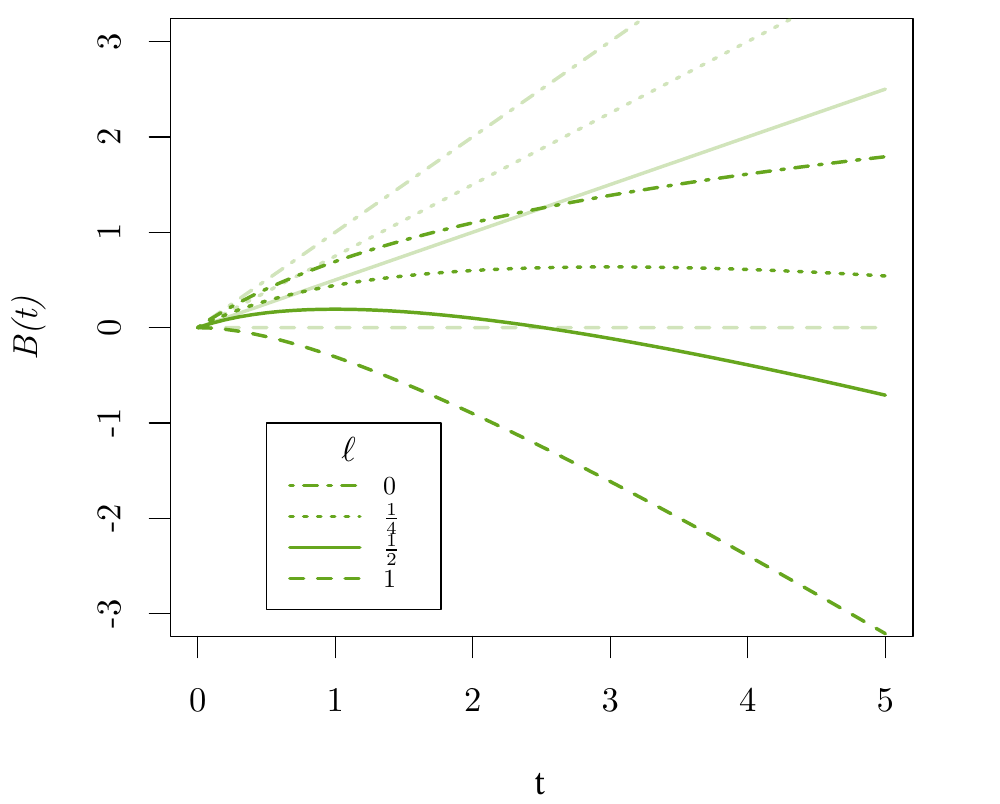}} 
		\caption{$\int_{0}^{t} \mathbb{E}\left[U_{1} \mid T^{1}\geq s\right] ds$ (thick lines), when $U_{1}+\ell \sim \Gamma(1,1)$  and $g(t)=t \mathbb{E}[U_{1}]$ (transparent lines). 
  }\label{CH6FigGamma}	
\end{figure}

Similar patterns can be observed for a continuous $U_{1}$ distribution, in which case the $\mathbb{E}\left[U_{1} \mid T^{1} \geq t\right]$ will keep decreasing. For example when $(U_{1}+\ell) \sim \Gamma(k, \theta)$, then $\mathcal{L}_{U_{1}+\ell}(c)= \frac{\theta k}{\theta c +1}$. By lemma \ref{CH6lemLP}
$$ \mathbb{E}\left[U_{1} \mid T^{1} \geq t\right]=\frac{\theta k}{\theta t +1} - \ell,  $$ 
and
$$ B(t) = k\log(\theta t + 1) - \ell t.$$ The $B(t)$ are presented over time in Figure \ref{CH6FigGamma} for $\theta=k=1$ and $\ell \in \{0, \tfrac{1}{4}, \tfrac{1}{2}, 1 \}$. 

Thus, if $U_{0} \independent U_{1}$, then the \mchd will be less or equal to the \chd due to the selection of $U_{1}$. Decreasing or constant $B(t)$ curves that at some point increase again can then not be explained by the selection of $U_{1}$ since individuals with less beneficial values of $U_{1}$ are expected to survive shorter. However, if $U_{0} \not \independent U_{1}$, such a pattern of the $B(t)$ curve can be due to the selection of $U_{1}$ as, by Theorem \ref{CH6th52}, $\mathbb{E}[f_{0}(t,U_{0})\mid T^{1} \geq t] \neq \mathbb{E}[f_{0}(t,U_{0})\mid T^{0} \geq t]$.  If $U_{0} \not \independent U_{1}$, the deviation of the \ohd from the \chd depends on $f_{0}$ and the joint distribution of $(U_{0},U_{1})$. 

\subsection{Dependent $U_{0}$ and $U_{1}$}\label{CH6sec:cop}
\begin{figure}[H]
	\centering
	\captionsetup{width=\textwidth}
	\includegraphics[width=\textwidth]{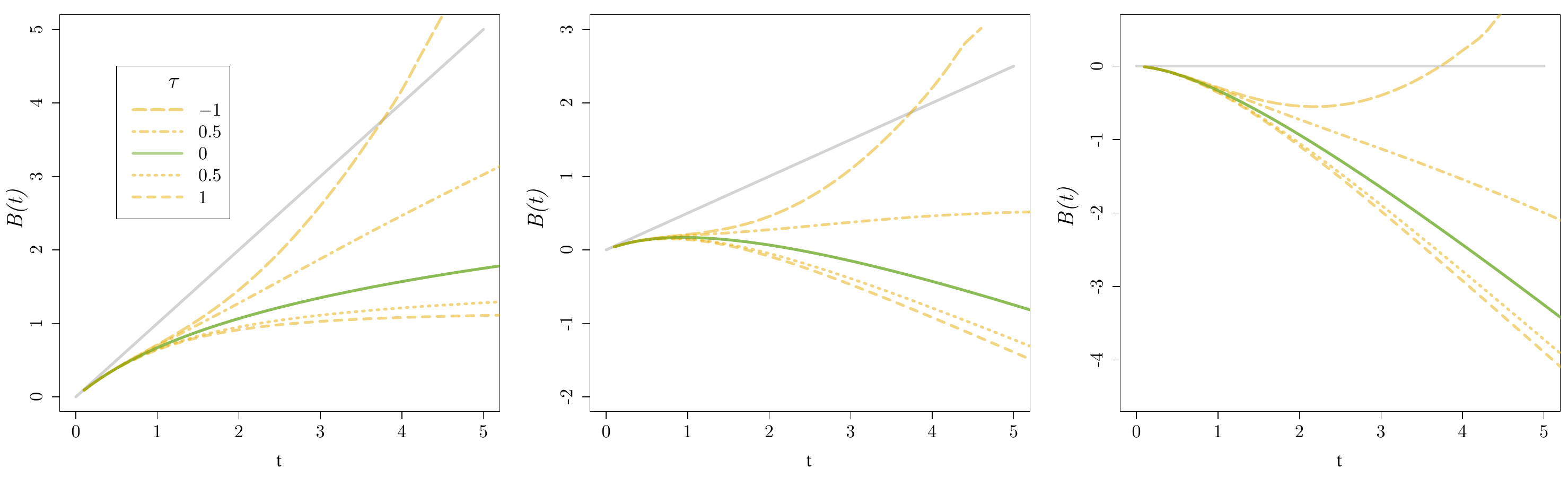} 
	\caption{Integrated hazard difference, $B(t)$, when $f_{0}(t,U_{0i})=\ell + U_{0i}t^{2}$, $U_{0i}\sim \Gamma(1,1)$, $U_{1i}+ \ell \sim \Gamma(1,1)$ for $\ell$ equal to $0$ (left), $\tfrac{1}{2}$ (middle) and $1$ (right) and different Kendall's $\tau$ correlation coefficients of $U_{0i}$ and $U_{1i}$. The lines for $\rho=0$ were already presented in Figure \ref{CH6FigGamma}. Furthermore, $g(t)=t \mathbb{E}[U_{1}]$ are presented as gray lines.  }\label{CH6figcop}	
\end{figure}

\begin{figure}[H]
	\centering
	\captionsetup{width=\textwidth}
	\includegraphics[width=\textwidth]{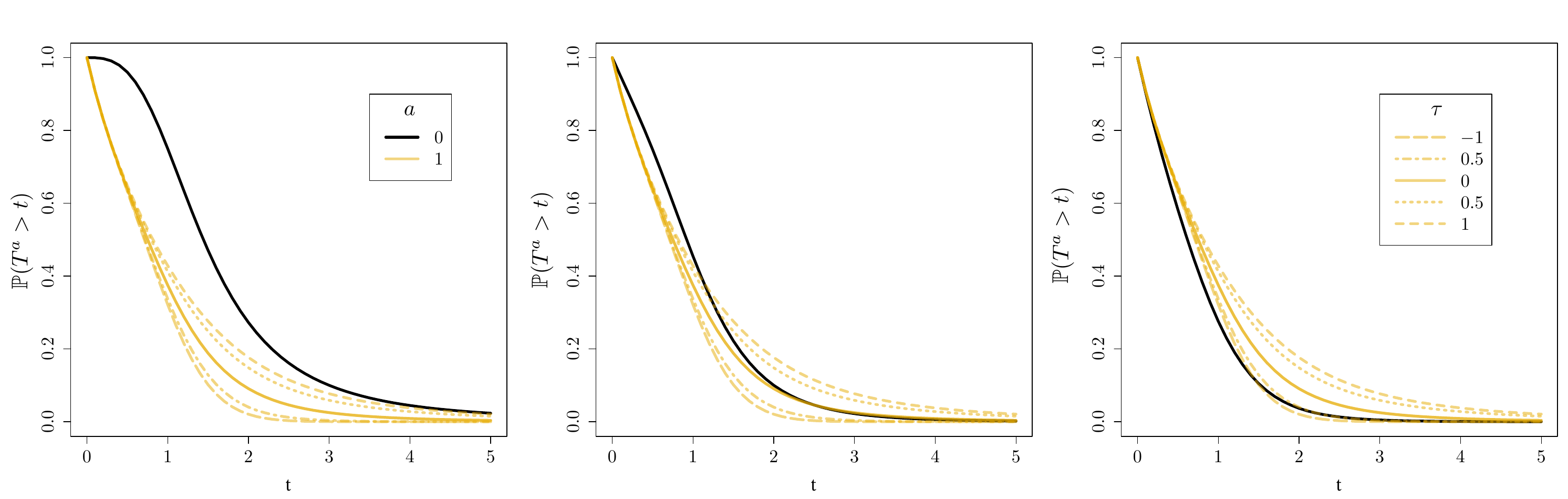} 
	\caption{	Survival curves for $Y^{1}$ and $Y^{0}$ when $f_{0}(t,U_{0i})=\ell + U_{0i}t^{2}$, $U_{0i}\sim \Gamma(1,1)$, $U_{1i}+ \ell \sim \Gamma(1,1)$ for $\ell$ equal to $0$ (left), $\tfrac{1}{2}$ (middle) and $1$ (right) and different Kendall's $\tau$ correlation coefficients of $U_{0i}$ and $U_{1i}$.  }\label{CH6figcopS}	
\end{figure}

The bivariate joint distribution function of $U_{0}$ and $U_{1}$, $F_{(U_{0},U_{1})}$, can be written using the marginal distribution functions and a copula $C$ \citep{Sklar1959}. As such, 
\begin{equation*}
F_{(U_{0},U_{1})}(u_{0},u_{1})=C\left(F_{U_{0}}(u_{0}), F_{U_{1}}(u_{1})\right)  
\end{equation*} and the Kendall's $\tau$ correlation coefficient of $U_{0}$ and $U_{1}$ can be written as a function of the Copula $C$ \citep{Nelsen2006}. To illustrate how the dependence can affect the \ohd for the settings presented in Figure \ref{CH6FigGamma}, we use a Gaussian copula \begin{equation*}
C(x,y)=\Phi_{2,\rho}(\Phi^{-1}(x), \Phi^{-1}(y)),
\end{equation*} where $\Phi$ and $\Phi_{2,\rho}$ are the standard normal and standard bivariate normal with correlation $\rho$ cumulative distribution functions, respectively. Then, 
$$F_{(U_{0},U_{1})}(u_{0},u_{1})= \Phi_{2,\rho}(\Phi^{-1}(F_{U_{0}}(u_{0}), \Phi^{-1}(F_{U_{1}}(u_{1})))).$$ For this example, we let 
$$f_{0}(t,U_{0i})=\ell + U_{0i}t^{2},$$ where $U_{0i} \sim \Gamma(1,1)$ so that the hazard is nonnegative for each individual. For \scalebox{0.9}{$\rho \in \left\{-1, \sin\left(-0.25\pi\right), 0, \sin\left(0.25\pi\right), 1\right \}$} (such that $\tau \in \{-1, -0.5, 0, 0.5, 1\}$) and $\ell \in \{0, \tfrac{1}{2}, 1 \}$, $\int_{0}^{t} \left( \mathbb{E}[ U_{1}\mid T^{1}\geq s] + \mathbb{E}[ U_{0}s^{2}\mid T^{1}\geq s] -\mathbb{E}[U_{0}s^{2} \mid T^{0}\geq s] \right) ds$ is presented in Figure \ref{CH6figcop}. The conditional expectations are derived empirically from simulations $(n=10,000)$, and the integral is approximated by taking discrete steps of size $0.1$. All programming codes used for this work can be found online at \url{https://github.com/RAJP93/CHD}. The survival curves of the potential outcomes can be found in Figure \ref{CH6figcopS}. 

\noindent The selection bias increases when $\rho>0$ (compared to $\rho=0$). On the other hand, for $\rho<0$,  the selection bias is smaller most of the time than for $\rho=0$ since favourable $U_{1}$ are expected to occur with unfavourable levels of $U_{0}$. Moreover, for $\rho=-1$, at larger $t$, we observe that the selection bias can even change sign, and the absolute value of the bias might be larger than for $\rho=0$ for those $t$. For $\rho \neq 0$, the \mchd might thus be larger than the \chd, so the \mchd is not a lower bound for the \chd. Note that if $U_{0} \not \independent U_{1}$, the integrated \mchd depends on the functional form of $f_{0}$. In Figure \ref{CH6figcopAlt} in Appendix \ref{CH6Btalt} the results for $f_{0}(t,U_{0i})=\ell + U_{0i}\frac{t^{2}}{20}$ are presented where the effect of the dependence is small.  

\section{Case study: the Radiation Therapy Oncology Group trial}
\begin{figure}[H]
	\centering
	\captionsetup{width=\textwidth}
	\includegraphics[width=0.7\textwidth]{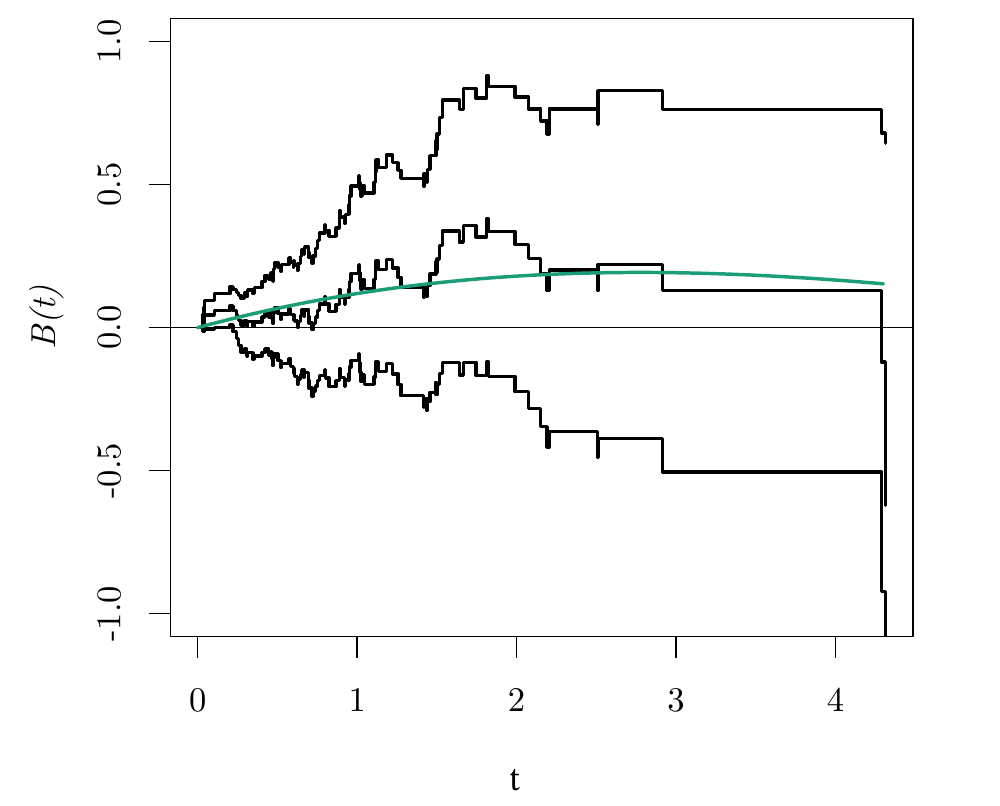} 
	\caption{Estimated $B(t)$ and corresponding $95\%$ confidence interval (black), including patients with tumours at one of the two sites considered. Furthermore, the expected evolution of $B(t)$ when 
		$ \lambda_{i}^{a}(t)=f_{0}(t,U_{0i})+U_{1i}$ and $U_{1i} \sim \text{BHN}(0.5, -0.1, 0.5, 0.4)$ is presented (green). }\label{CH6FigCS1}	
\end{figure}

\begin{figure}[H]
	\centering
	\captionsetup{width=\textwidth}
	\includegraphics[width=\textwidth]{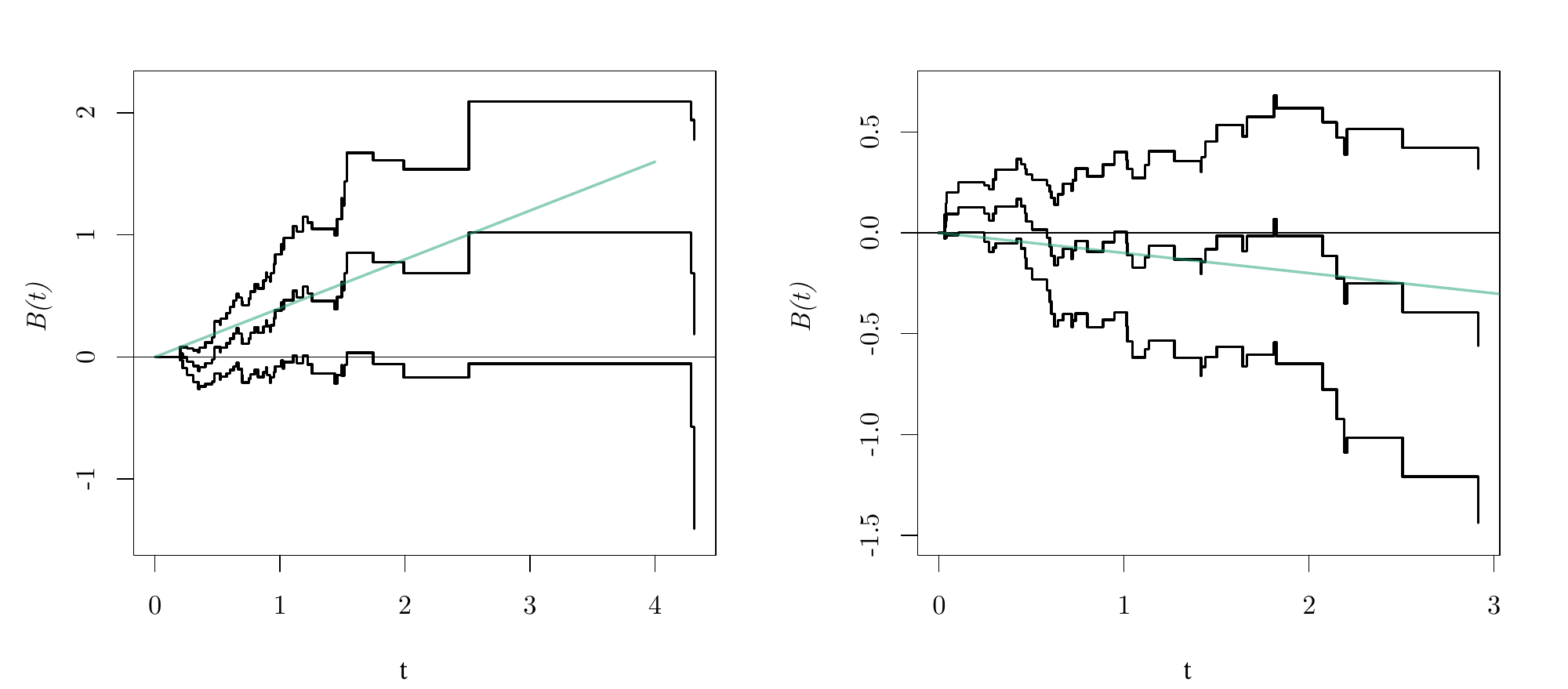} 
	\caption{Estimated $B(t)$ and corresponding $95\%$ confidence interval (black) for patients with tumours located at the faucial arch (left) and pharyngeal tongue (right), respectively. Furthermore, the $B(t)$ for a homogeneous population is presented (green) equal to $0.4t$ (right) and $-0.1t$ (left) for comparison }\label{CH6FigCS2}	
\end{figure}

Next, we consider a large clinical trial carried out by the Radiation Therapy Oncology Group as described by \citet[Section 1.1.2 and Appendix A]{Kalbfleisch2002} and also presented by \citet{Aalen1989}. From the patients with squamous cell carcinoma (a form of skin cancer) of 15 sites in the mouth and throat from 16 participating institutions, our focus is only on two sites (faucial arch and pharyngeal tongue) and patients from the six largest institutions. All participants were randomly assigned to radiation therapy alone or combined with a chemotherapeutic agent. So, we are interested in the causal effect of the chemotherapeutic agent in addition to radiation therapy on survival. If the causal mechanism can be parameterized with SCM \eqref{CH6SCMsurv} without effect heterogeneity, i.e.~ 
$$ \lambda_{i}^{a}(t)=f_{0}(t,U_{0i})+f_{1}(t,a),$$ and the randomization was properly executed, implying $N_{Ai} \independent U_{0i}, N_{Ti}$, then, by Theorem \ref{CH6th52}, the \ohd equals the \chd. Furthermore, the \chd can be unbiasedly estimated by fitting Aalen's additive hazard model. We did so by using the \texttt{aalen()} function from the package \texttt{timereg} in \texttt{R}. The estimated cumulative regression function (and a corresponding $95\%$ confidence interval) of treatment combined with a chemotherapeutic agent is presented by the black lines in Figure \ref{CH6FigCS1}. 
\noindent In the absence of effect heterogeneity (ignoring the statistical uncertainty), one could now conclude that initially adding the chemotherapy is expected to harm a patient as $B(t)$ takes on positive values and that the exposure loses its harmful effect over time as the derivative of $B(t)$ decreases over time. Following a similar reasoning, \citet{Aalen2008} concluded for the effect of N-stage (an index of lymph node metastasis) on survival based on the same dataset (while also including patients with a tumour located at the tonsillar fossa) that \textit{``The regression plot shows
	that this} [non-significant P-value for a zero-effect of N-stage from a Cox analysis] \textit{is due to a strong initial positive effect being "watered down" by a lack of, or even a slightly negative effect after one year. Hence, not taking into consideration the change in effect over time may lead to missing significant effects."}. 
However, when in reality
$$ \lambda_{i}^{a}(t)=f_{0}(t,U_{0i})+f_{1}(t,U_{1i},a),$$
the observed time-varying effect can also result from the modifier $U_{1i}$ selection. For example, when 
$ \lambda_{i}^{a}(t)=f_{0}(t,U_{0i})+U_{1i}$, and $U_{1i} \sim \text{BHN}(0.5, -0.1, 0.5, 0.4)$, by Theorem \ref{CH6th52}, this pattern is expected (see the green line in Figure \ref{CH6FigCS1}) while the causal effect is constant for each individual. The \chd equals $0.15$ at each time point, but over time individuals that are harmed by the chemotherapy ($U_{1i}=0.4$) are less likely to survive so that the \mchd converges towards $-0.1$ (the effect for individuals that benefit from the chemotherapy). When we perform a stratified analysis by site in the oropharynx (where randomization remains), we observe that the effect of chemotherapy might have opposite effects for tumours located in the faucial arch and on the pharyngeal tongue, see Figure \ref{CH6FigCS2}. The tumour location could thus be the individual modifier underlying the BHN distribution. For this case study, we cannot be sure whether the effect of chemotherapy depends on the tumour location due to statistical uncertainty. However, it became clear that when statistical uncertainty is not the issue, it will be impossible to distinguish between a time-varying causal effect and a selection effect (of an unmeasured modifier) from data. Both phenomena can give rise to the same \mchd curves.  

\section{Discussion and concluding remarks}
The additive hazard model gives better interpretable estimates of causal effects than the proportional hazard model \citep{Aalen2015}. As discussed by \citet{Aalen2015}, the model assumes that the additive part of the hazard involving the exposure (or treatment) is not affected by any other individual feature. Otherwise, if such effect heterogeneity at the hazard scale exists, we have shown that the \ohd is a biased measure for the \chd of interest. A time-varying \ohd can thus be the result of either an actual time-varying causal effect or of the selection of favourable effect-modifier levels over time. Therefore, it is impossible to distinguish these scenarios based on data without making untestable assumptions. It is important to remark that for cause-effect relations that can be parameterized with SCM \eqref{CH6SCMsurv} where $U_{1}$ is degenerate (in which case the \ohd can be an unbiased estimator of the \chd), contrary to the individual hazard differences, the difference of the potential survival times, $T^{1}-T^{0}$ can be heterogeneous. So, heterogeneous effects can also exist under Aalen's additive hazard model. 

In the case study, we have illustrated that one should be very careful when concluding that the effect decreases over time based on the cumulative regression function, as this might result from the selection. The size of the bias depends on how much  the distribution $F_{U_{1}\mid T^{1} \geq t}$ changes over time. In case the $U_{1}$ is low in variability, the size of the bias will be small. When analyzing data from an RCT with an additive hazard model, it can thus be helpful to adjust for potential effect modifiers to reduce the remaining variability of unmeasured effect modifiers. 

Even in the absence of confounding, the hazard difference and the hazard ratio (as discussed in our companion paper \citet{Post2022b}) are biased estimators of the causal effect. Instead, contrasts of the survival probabilities, the median, or the restricted mean survival time, can be unbiased estimators of causal estimands and should thus be used to quantify the causal effect on time-to-event outcomes as suggested by others \citep{Hernan2010, Bartlett2020, Stensrud2019, Young2020}.

		\bibliographystyle{spbasic}
		\bibliography{CHR}

\newpage
\appendix

\section{Proofs}\label{CH6CH6:proofs}
\subsection{Proof of Lemma \ref{CH6l5.1}}
\begin{proof}
	\begin{align*}
	&\mathbb{E}\left[f_{0}(t,U_{0})\mid T^{a} \geq t\right]\\ &=\int f_{0}(t,U_{0}) f_{U_{0}\mid T^{a} \geq t}(u_{0}) du_{0}\\
	&=\int f_{0}(t,U_{0}) \frac{\mathbb{P}(T^a \geq t\mid U_{0}=u_{0})}{\mathbb{P}(T^{a}\geq t)}f_{U_{0}}(u_{0}) du_{0}\\
	&  = \int f_{0}(t,U_{0}) \frac{\int\exp(-(\int_{0}^{t}f_{0}(u_{0},s)ds+\int_{0}^{t}f_{1}(u_{1},a,s)ds ))f_{U_{1}|U_{0}=u_{0}}(u_{1})du_{1}}{\int\int \exp(-(\int_{0}^{t}f_{0}(k_{0},s)ds + \int_{0}^{t}f_{1}(k_{1},a,s)ds))f_{U_{1}\mid U_{0}=k_{0}}(k_{1})dk_{1}f_{U_{0}}(k_{0}) dk_{0}} f_{U_{0}}(u_{0}) du_{0}\\
	&  = \int f_{0}(t,U_{0}) \frac{\exp(-\int_{0}^{t}f_{0}(u_{0},s)ds) \left(\int\exp(-\int_{0}^{t}f_{1}(u_{1},a,s)ds )f_{U_{1}|U_{0}=u_{0}}(u_{1})du_{1}\right)}{\int\exp(-\int_{0}^{t}f_{0}(k_{0},s)ds ) \left(\int \exp(- \int_{0}^{t}f_{1}(k_{1},a,s)ds)f_{U_{1}\mid U_{0}=k_{0}}(k_{1})dk_{1}\right)f_{U_{0}}(k_{0}) dk_{0}} f_{U_{0}}(u_{0}) du_{0}.
	\end{align*} 	Moreover, if $U_{0}\independent U_{1}$, then
	\begin{align*}
	&\mathbb{E}\left[f_{0}(t,U_{0})\mid T^{a} \geq t\right]\\
	&  = \int f_{0}(t,U_{0}) \frac{\exp(-\int_{0}^{t}f_{0}(u_{0},s)ds) \left(\int\exp(-\int_{0}^{t}f_{1}(u_{1},a,s)ds )f_{U_{1}}(u_{1})du_{1}\right)}{\int\exp(-\int_{0}^{t}f_{0}(k_{0},s)ds ) \left(\int \exp(- \int_{0}^{t}f_{1}(k_{1},a,s)ds)f_{U_{1} }(k_{1})dk_{1}\right)f_{U_{0}}(k_{0}) dk_{0}} f_{U_{0}}(u_{0}) du_{0}\\
	&  = \int f_{0}(t,U_{0}) \frac{\exp(-\int_{0}^{t}f_{0}(u_{0},s)ds) \left(\int\exp(-\int_{0}^{t}f_{1}(u_{1},a,s)ds )f_{U_{1}}(u_{1})du_{1}\right)}{\left(\int\exp(-\int_{0}^{t}f_{0}(k_{0},s)ds) f_{U_{0}}(k_{0}) dk_{0}\right) \left(\int \exp(- \int_{0}^{t}f_{1}(k_{1},a,s)ds)f_{U_{1} }(k_{1})dk_{1}\right)} f_{U_{0}}(u_{0}) du_{0}\\
	&= \int f_{0}(t,U_{0}) \frac{\exp(-\int_{0}^{t}f_{0}(u_{0},s)ds) }{\int\exp(-\int_{0}^{t}f_{0}(k_{0},s)ds) f_{U_{0}}(k_{0}) dk_{0}} f_{U_{0}}(u_{0}) du_{0}\\
	&= \mathbb{E}\left[f_{0}(t,U_{0})\mid T^{0} \geq t\right].  
	\end{align*} 
\end{proof}

\subsection{Proof of Theorem \ref{CH6th52}}
\begin{proof}
	
	By causal consistency \citep{Hernan2019},
	$$
	\lim_{h\rightarrow 0}h^{-1}\mathbb{P}\left(T \in [t,t+h) \mid T \geq t, A=a \right) =  \lim_{h\rightarrow 0}h^{-1}\mathbb{P}\left(T^{a} \in [t,t+h) \mid T^{a} \geq t, A=a \right)
	$$
	As $N_{A} \independent U_{0}, U_{1}, N_{T}$, there is no confounding and $T^{a} \independent A$, so that 
	
	$$
	\lim_{h\rightarrow 0}h^{-1}\mathbb{P}\left(T \in [t,t+h) \mid T \geq t, A=a \right) =  \lim_{h\rightarrow 0}h^{-1}\mathbb{P}\left(T^{a} \in [t,t+h) \mid T^{a} \geq t, \right).
	$$ By the law of total probability, 
	$$
\lim_{h\rightarrow 0}h^{-1}\mathbb{P}\left(T^{a} \in [t,t+h) \mid T^{a} \geq t \right) =
	\int\lim_{h\rightarrow 0}h^{-1}\mathbb{P}\left(T^{a} \in [t,t+h) \mid T^{a} \geq t, U_{0}, U_{1} \right)dF_{U_{0}, U_{1}\mid T^{a} \geq t}$$
	First, we focus on the integrand,  \begin{align*}
	&h^{-1}\mathbb{P}\left(T^{a} \in [t,t+h) \mid T^{a} \geq t, U_{0}, U_{1} \right)\\  &=h^{-1}\frac{\mathbb{P}\left(T^{a} \geq  t \mid  U_{0}, U_{1} \right) - \mathbb{P}\left(T^{a} \geq  t + h \mid  U_{0}, U_{1} \right)}{\mathbb{P}\left(T^{a} \geq  t \mid  U_{0}, U_{1} \right)} \\ 
	&=  h^{-1}\left(1-\frac{\mathbb{P}\left(T^{a} \geq  t + h \mid  U_{0}, U_{1} \right)}{\mathbb{P}\left(T^{a} \geq  t \mid  U_{0}, U_{1} \right)}\right) \\ 
	&=  h^{-1}\left(1-\frac{\exp\left(- \int_{0}^{t+h} f_{0}(U_{0},s)+f_{1}(U_{1},a,s) ds \right)}{\exp\left(- \int_{0}^{t} f_{0}(U_{0},s)+f_{1}(U_{1},a,s) ds \right)}\right) \\
	&=  h^{-1}\left(1-\exp\left(- \int_{t}^{t+h} f_{0}(U_{0},s)+f_{1}(U_{1},a,s) ds \right)\right) 
	\end{align*} For monotonic conditional hazard functions, if $h_{2}<h_{1}$, then \begin{equation*}
	 h_{1}^{-1}\left(1- \exp\left(- \int_{t}^{t+h_{1}} f_{0}(U_{0},s)+f_{1}(U_{1},a,s) ds \right)\right) \leq  h_{2}^{-1}\left(1- \exp\left(- \int_{t}^{t+h_{2}} f_{0}(U_{0},s)+f_{1}(U_{1},a,s) ds \right)\right)
	\end{equation*} or 
	\begin{equation*}
	 h_{1}^{-1}\left(1- \exp\left(- \int_{t}^{t+h_{1}} f_{0}(U_{0},s)+f_{1}(U_{1},a,s) ds \right)\right) \geq  h_{2}^{-1}\left(1- \exp\left(- \int_{t}^{t+h_{2}} f_{0}(U_{0},s)+f_{1}(U_{1},a,s) ds \right)\right)
	\end{equation*} as the average integrated conditional hazard over the interval increases (or decreases). Then, the limit and integral can be interchanged by directly applying the monotone convergence theorem.
	
	For non-monotone conditional hazard functions, let $h \leq \tilde{{}h}$, and note 
	\begin{equation*}h^{-1}\mathbb{P}\left(T^{a} \in [t,t+h) \mid T^{a} \geq t, U_{0}, U_{1} \right) \leq h^{-1}\left(1- \exp\left(- h\left( f_{0}(t,U_{0}^{*})+f_{1}(U_{1},a,t^{*})\right)\right)\right),\end{equation*} \noindent where
	\begin{equation*}
	\max_{s \in (t, t + \tilde{{}h})} f_{0}(U_{0},s)+f_{1}(U_{1},a,s) = f_{0}(t,U_{0}^{*})+f_{1}(U_{1},a,t^{*}).
	\end{equation*} Using the power series definition of the exponential function, \begin{align*}
	&h^{-1}\mathbb{P}\left(T^{a} \in [t,t+h) \mid T^{a} \geq t, U_{0}, U_{1} \right)\\
 &\leq h^{-1}\left(1- \frac{1}{\sum_{k=0}^{\infty}h^{k}  (f_{0}(t,U_{0}^{*})+f_{1}(U_{1},a,t^{*}))^{k}\tfrac{1}{k!}}\right)\\
	&= h^{-1} \frac{\sum_{k=1}^{\infty}h^{k}  (f_{0}(t,U_{0}^{*})+f_{1}(U_{1},a,t^{*}))^{k}\tfrac{1}{k!}}{\sum_{k=0}^{\infty}h^{k}  (f_{0}(t,U_{0}^{*})+f_{1}(U_{1},a,t^{*}))^{k}\tfrac{1}{k!}}\\
	&= f_{0}(t,U_{0}^{*})+f_{1}(U_{1},a,t^{*}) \frac{\sum_{k=1}^{\infty}h^{k-1}  (f_{0}(t,U_{0}^{*})+f_{1}(U_{1},a,t^{*}))^{k-1}\tfrac{1}{k!}}{\sum_{k=0}^{\infty}h^{k}  (f_{0}(t,U_{0}^{*})+f_{1}(U_{1},a,t^{*}))^{k}\tfrac{1}{k!}}\\
	&= f_{0}(t,U_{0}^{*})+f_{1}(U_{1},a,t^{*}) \frac{\sum_{k=0}^{\infty}h^{k}  (f_{0}(t,U_{0}^{*})+f_{1}(U_{1},a,t^{*}))^{k}\tfrac{1}{(k+1)!}}{\sum_{k=0}^{\infty}h^{k}  (f_{0}(t,U_{0}^{*})+f_{1}(U_{1},a,t^{*}))^{k}\tfrac{1}{k!}}\\
	&< f_{0}(t,U_{0}^{*})+f_{1}(U_{1},a,t^{*}).
	\end{align*}  
	Moreover, $\mathbb{E}\left[f_{0}(t,U_{0}^{*})+f_{1}(U_{1},a,t^{*})\mid T^{a} \geq t\right]<\infty$ when   $\mathbb{E}[f_{0}(t,U_{0}+h)+f_{1}(U_{1},a,t+h)\mid T^{a} \geq t]<\infty$ for all $h \in (0,\tilde{{}h})$. By application of the dominated convergence theorem, we can change the order of the limit and integral and conclude, 
	\begin{equation*}
	\lim_{h\rightarrow 0}h^{-1}\mathbb{P}\left(T^{a} \in [t,t+h) \mid T^{a} \geq t \right)
	= \mathbb{E}\left[f_{0}(t,U_{0})+f_{1}(U_{1},a,t)\mid T^{a} \geq t\right]
	\end{equation*}  
	As $U_{0} \independent U_{1}$, by Lemma \ref{CH6l5.1}, 
	$$
	\mathbb{E}[f_{0}(t,U_{0}) + f_{1}(t,U_{1},1)\mid T^{1} \geq t]-\mathbb{E}[f_{0}(t,U_{0}) \mid T^{0} \geq t]=\mathbb{E}\left[f_{1}(t,U_{1},1)\mid T^{1} \geq t\right],
	$$ so 
	\begin{align*}
	&\lim_{h\rightarrow 0}h^{-1}\mathbb{P}\left(T \in [t,t+h) \mid T \geq t, A=1 \right)-\lim_{h\rightarrow 0}h^{-1}\mathbb{P}\left(T\in [t,t+h) \mid T \geq t, A=0 \right)\\
	& = \mathbb{E}\left[f_{1}(t,U_{1},1)\mid T^{1} \geq t\right]. 
	\end{align*}
\end{proof}

\subsection{Proof of Lemma \ref{CH6lemLP}}
\begin{proof}
	By Bayes rule, the probability density $f_{U_{1}\mid T^{1} \geq t}$ equals
	\begin{align*}
	&f_{U_{1}\mid T^{1} \geq t}(u_{1})\\ 
        &= \frac{\mathbb{P}(T^{1}\geq t\mid U_{1}=u_{1})f(u_{1})}{\int \mathbb{P}(T^{1}\geq t\mid U_{1})dF_{U_{1}}}\\
	&= \frac{\int\exp(-\left(\int_{0}^{t}f_{0}(U_{0},s)+u_{1}f_{1}(a,s) ds)\right)f_{U_{1}}(u_{1})dF_{U_{0}}}{\int\int \exp(-\left(\int_{0}^{t}f_{0}(U_{0},s)+U_{1}f_{1}(a,s)ds\right))dF_{U_{1}}dF_{U_{0}}}\\
	&= \frac{\exp(-\int_{0}^{t}u_{1}f_{1}(a,s) ds)f_{U_{1}}(u_{1})\int\exp(-\int_{0}^{t}f_{0}(U_{0},s) ds)dF_{U_{0}}}{\int\exp(-U_{1}f_{1}(a,s)ds)dF_{U_{1}}\int \exp(-\int_{0}^{t}f_{0}(U_{0},s)ds)dF_{U_{0}}}\\
	&= \frac{\exp(-\int_{0}^{t}u_{1}f_{1}(a,s) ds)f_{U_{1}}(u_{1})}{\int\exp(-U_{1}f_{1}(a,s)ds)dF_{U_{1}}}.
	\end{align*} \noindent So that the Laplace transform of $f_{U_{1}\mid T^{1} \geq t}$ can be written as
	\begin{align*}
	\mathcal{L}_{U_{1}\mid T^{1} \geq t}(c) &= \mathbb{E}[\exp(-U_{1}c)\mid T^{1} \geq t]\\
	&= \int \exp(-u_{1}c)  f_{U_{1}\mid T^{1} \geq t}(u_{1}) du_{1}\\
	&= \int \exp(-u_{1}c) \frac{\exp(-\int_{0}^{t}u_{1}f_{1}(1,s) ds)f_{U_{1}}(u_{1})}{\int\exp(-U_{1}f_{1}(1,s)ds)dF_{U_{1}}}du_{1}\\
	&= \int \frac{\exp(-u_{1}(c+\int_{0}^{t}f_{1}(1,s) ds))f_{U_{1}}(u_{1})}{\int\exp(-U_{1}f_{1}(1,s)ds)dF_{U_{1}}}du_{1}\\
	&=\int \exp(-u_{0}c) \frac{\mathbb{E}\left[\exp(-\int_{0}^{t}u_{1}f_{1}(1,s) ds)\right]}{\mathbb{E}\left[\exp(-U_{1}f_{1}(1,s)ds)\right]}du_{1}\\
	&= \frac{\mathcal{L}_{U_{1}}(c+\int_{0}^{t}f_{1}(1,s) ds)}{\mathcal{L}_{U_{1}}(\int_{0}^{t}f_{1}(1,s) ds)}. 
	\end{align*}
\noindent Since for a random variable $X$, $\mathbb{E}[X]= -\mathcal{L}_{X}^{'}(0)$, 
	\begin{equation*}
	\mathbb{E}[U_{1}\mid T^{1} \geq t] = -\mathcal{L}_{U_{1}\mid T^{1} \geq t}^{'}(0) = -\frac{\mathcal{L}_{U_{1}}^{'}(\int_{0}^{t}f_{1}(1,s)ds)}{\mathcal{L}_{U_{1}}(\int_{0}^{t}f_{1}(1,s)ds)}. 
	\end{equation*}
\end{proof}

\newpage
\section{Figures dependent $U_{0}$ and $U_{1}$}\label{CH6Btalt}
\begin{figure}[H]
	\centering
	\captionsetup{width=\textwidth}
	\includegraphics[width=\textwidth]{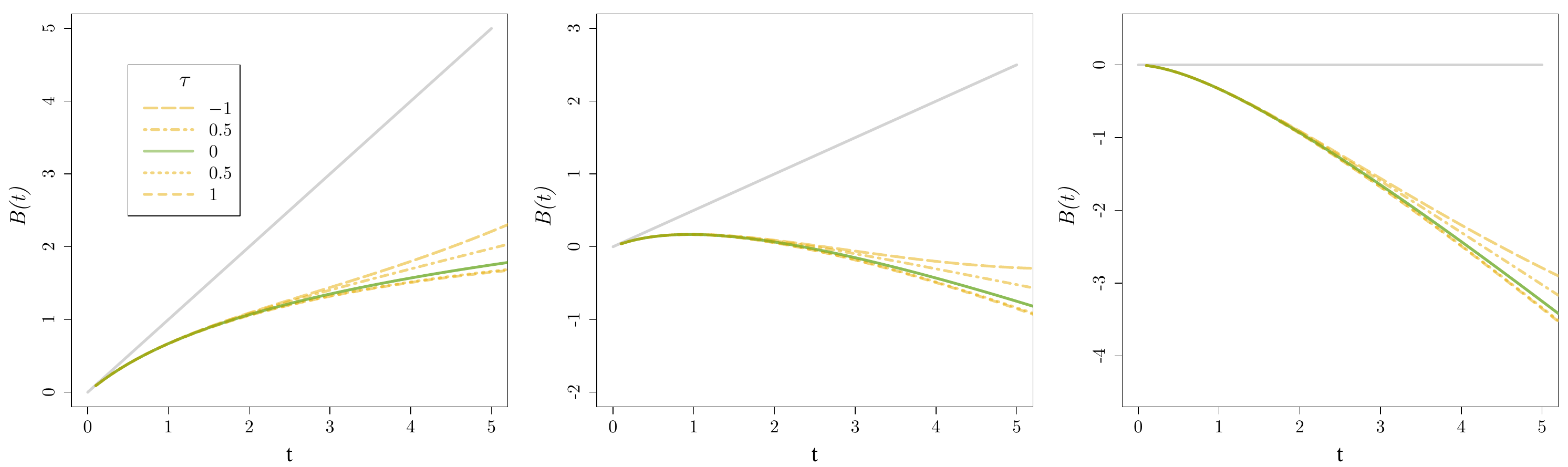} 
	\caption{	Integrated hazard difference, $B(t)$, when $f_{0}(t,U_{0i})=\ell + U_{0i}\tfrac{t^{2}}{20}$, $U_{0i}\sim \Gamma(1,1)$, $U_{1i}+ \ell \sim \Gamma(1,1)$ for $\ell$ equal to $0$ (left), $\tfrac{1}{2}$ (middle) and $1$ (right) and different Kendall's $\tau$ correlation coefficients of $U_{0i}$ and $U_{1i}$. The lines for $\rho=0$ were already presented in Figure \ref{CH6FigGamma}. Furthermore, $g(t)=t \mathbb{E}[U_{1}]$ are presented as gray lines.  }\label{CH6figcopAlt}	
\end{figure}

The additional figures for $f_{0}(t,U_{0i})=\ell + U_{0}\frac{t^{2}}{20}$. In Figure \ref{CH6figcopAlt}, the lines for $\rho=0.5$ and $\rho=1$ do overlap. 

\begin{figure}[H]
	\centering
	\captionsetup{width=\textwidth}
	\includegraphics[width=\textwidth]{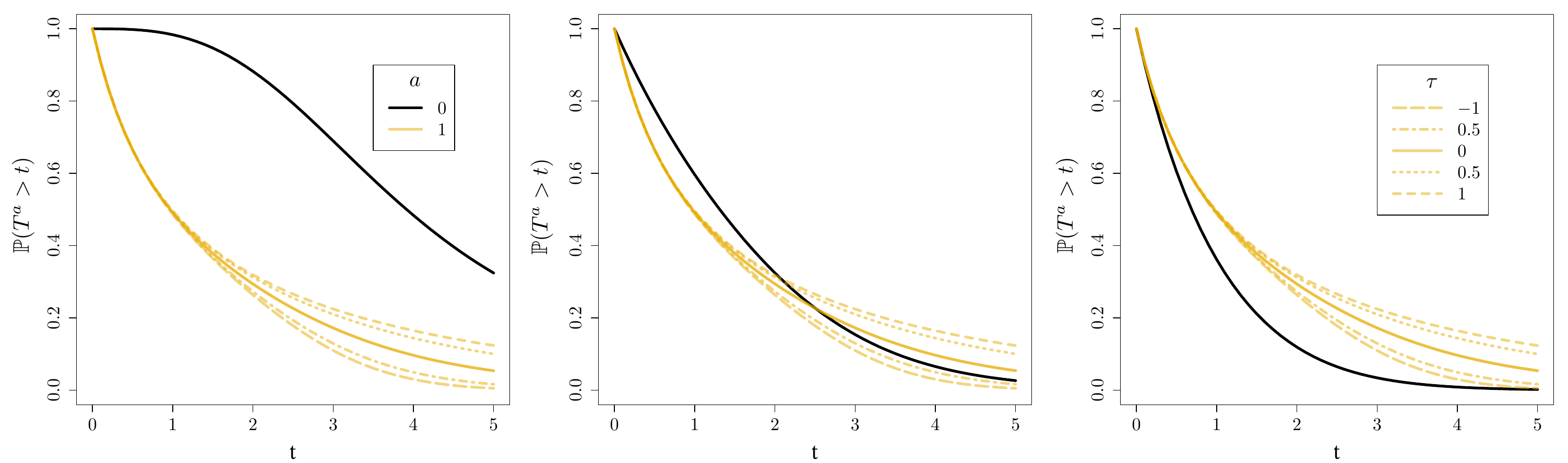} 
	\caption{		Survival curves for $Y^{1}$ and $Y^{0}$ when $f_{0}(t,U_{0i})=\ell + U_{0i}\tfrac{t^{2}}{20}$, $U_{0i}\sim \Gamma(1,1)$, $U_{1i}+ \ell \sim \Gamma(1,1)$ for $\ell$ equal to $0$ (left), $\tfrac{1}{2}$ (middle) and $1$ (right) and different Kendall's $\tau$ correlation coefficients of $U_{0i}$ and $U_{1i}$.   }\label{CH6figcopAltS}	
\end{figure}

\end{document}